\documentclass[10pt]{article}
\usepackage{latexsym}
\usepackage{amssymb}
\usepackage{amsthm,amsmath}
\usepackage{graphicx}

\newtheorem{lemma}{Lemma}
\newtheorem{theorem}[lemma]{Theorem}
\newcommand{\R}{\mathbb R}
\newcommand{\E}{\mathbb E}
\newcommand{\PP}{\mathbb P}

\title{A short proof of Paouris' inequality}

\date{}
\author{
Rados{\l}aw Adamczak
 \thanks{Research partially supported by MNiSW Grant no. N N201 397437.}
 \and  Rafa{\l} Lata{\l}a
  \thanks{Research partially supported by MNiSW Grant no. N N201 397437.}
 \and  Alexander E. Litvak
\thanks{Research partially supported by  the
E.W.R. Steacie Memorial Fellowship.}
\and Krzysztof Oleszkiewicz
  \thanks{Research partially supported by MNiSW Grant no. N N201 397437.}
\and Alain Pajor
\thanks{Research partially supported by the ANR project
ANR-08-BLAN-0311-01.}
\and  Nicole  Tomczak-Jaegermann
\thanks{This author holds the Canada Research Chair in
  Geometric Analysis.}
}

\newcommand\address{\noindent\leavevmode%
\noindent
Rados{\l}aw  Adamczak, \\
Institute of Mathematics, \\
University of Warsaw, \\
Banacha 2, 02-097 Warszawa, Poland\\
 \texttt{\small
e-mail:  radamcz@mimuw.edu.pl}

\medskip
\noindent
Rafa{\l} Lata{\l}a, \\
Institute of Mathematics, \\
University of Warsaw, \\
Banacha 2, 02-097 Warszawa, Poland\\
 \texttt{\small
e-mail:   rlatala@mimuw.edu.pl}

\medskip
\noindent
Alexander E. Litvak, \\
Dept.~of Math.~and Stat.~Sciences,\\
University of Alberta, \\
Edmonton, Alberta, Canada, T6G 2G1.\\
\texttt{\small%
e-mail:   alexandr@math.ualberta.ca}

\medskip
\noindent
Krzysztof Oleszkiewicz, \\
Institute of Mathematics, \\
University of Warsaw, \\
Banacha 2, 02-097 Warszawa, Poland\\
 \texttt{\small
e-mail: koles@mimuw.edu.pl}

\medskip
\noindent
Alain  Pajor, \\
Universit\'{e} Paris-Est\\
\'{E}quipe d'Analyse et Math\'{e}matiques Appliqu\'ees, \\
5, boulevard Descartes,
Champs sur Marne,\\
77454 Marne-la-Vall\'{e}e,  Cedex 2, France\\
\texttt{\small%
e-mail: Alain.Pajor@univ-mlv.fr }

\medskip
\noindent
Nicole  Tomczak-Jaegermann, \\
Dept.~of Math.~and Stat.~Sciences,\\
University of Alberta, \\
Edmonton, Alberta, Canada, T6G 2G1.\\
\texttt{\small%
e-mail:    nicole.tomczak@ualberta.ca}
}

\begin{document}
\maketitle

\begin{abstract}
We give a short proof of a result of G.~Paouris on
the tail  behaviour of the Euclidean norm $|X|$ of  an isotropic
log-concave random vector $X\in\R^n,$
stating that for every $t\geq 1$,
\[\PP \big( |X|\geq ct\sqrt n\big)\leq \exp(-t\sqrt n).\]
More precisely we show that  for any log-concave random vector $X$
and any $p\geq 1$,
\[(\E|X|^p)^{1/p}\sim \E |X|+\sup_{z\in S^{n-1}}(\E |\langle
z,X\rangle|^p)^{1/p}.\]
\end{abstract}

(Secondary)   \\


\newpage{}

\section{Introduction}

Let $X$ be a random vector in the Euclidean space $\R^n$ equipped with
its Euclidean norm $|\,\cdot\,|$ and its scalar product
$\langle\cdot,\cdot\rangle$.
Assume that $X$ has a log-concave distribution (a typical
example of such a distribution is a random vector uniformly
distributed on a convex body). Assume further that it is centered
and its covariance matrix is the identity -- such a random
vector will be called {\em isotropic}. A famous and important
result of (\cite{Pa}, Theorem~1.1) states that

\begin{theorem}
\label{theorem:paouris}
There exists an absolute constant $c>0$ such that if $X$ is an
isotropic log-concave random vector in $\R^n$, then for every $t\geq 1$,
\[\PP \big( |X|\geq ct\sqrt n\big)\leq \exp(-t\sqrt n).\]
\end{theorem}

This result had a huge impact on the study of log-concave measures
and has a lot of applications in that subject.

A Borel  probability
measure on $\R^n$ is called log-concave if for all
$0<\theta<1$ and  all
compact sets $A,B\subset\R^n$ one has
\[
\mu( (1-\theta) A+\theta B)\geq \mu(A)^{1-\theta}\mu(B)^\theta.
\]
We refer to \cite{Bo1,Bo2} for a general study of this class of
measures. Clearly, the affine image of a log-concave probability is
also log-concave.  The Euclidean norm of an $n$-dimensional
log-concave random vector has moments of all orders (see
\cite{Bo1}). A log-concave probability is supported on some convex
subset of an affine subspace where it has a density.  In particular
when the support of the probability generates the whole space $\R^n$
(in which case we talk, in short, about full-dimensional probability)
a characterization of Borell (see \cite{Bo1, Bo2}) states that the
probability is absolutely continuous with respect to the Lebesgue
measure and has a density which is log-concave. We say that a random
vector is log-concave if its distribution is a log-concave measure.

Let $X\in\R^n$, be a random vector, denote the weak $p$-th moment of
$X$ by
\[
\sigma_{p}(X)=\sup_{z\in S^{n-1}}(\E |\langle z,X\rangle|^p)^{1/p}.
\]
The purpose of this article is to give a short proof of the following
\begin{theorem}
\label{theorem:imprPaouris}
For any log-concave random vector $X\in\R^n$ and any $p \ge 1$,
\[
(\E|X|^p)^{1/p}\leq C\big(\E |X|+\sigma_p(X)\big),
\]
where $C$ is an absolute positive constant.
\end{theorem}

This result may be deduced directly from Paouris' work in \cite{Pa}.
Indeed, it is a consequence of Theorem~8.2 combined with Lemma~3.9 in
that paper.  As formulated here, Theorem \ref{theorem:imprPaouris} first
appeared in \cite{cras_allpt} (Theorem 2).  Note that because
trivially a converse inequality is valid (with constant $1/2$), Theorem
\ref{theorem:imprPaouris} states in fact an equivalence for
$(\E|X|^p)^{1/p}$.

It is noteworthy  that
the following strengthening
of Theorem~\ref{theorem:imprPaouris} is still open:
$(\E|X|^p)^{1/p}\leq \E|X|+C\sigma_p(X)$,
where $C$ is an absolute positive constant.

\smallskip

If $X$ is a log-concave random vector, then so is $\langle z,X\rangle$
for every $z\in S^{n-1}$.  It follows that there exists an absolute
constant $C' > 0$ such that for any $p\geq 1$, $\sigma_{p}(X)\leq C'
p\,\sigma_{2}(X)$ (\cite{Bo1}).  (In fact one can deduce this
inequality with $C'= 1$ from \cite{bmp} or from Remark~5 in
\cite{klo}; see also Remark~1 following Theorem~3.1 in \cite{ALLPT}.)
%
If moreover $X$ is  isotropic, then $\E |X|\leq (\E |X|^2)^{1/2}=\sqrt n$
and $\sigma_2(X)= 1$;  thus
\[
(\E|X|^{p})^{1/p}\le C\big(\sqrt n+C'p\big).
\]
{}From Markov's inequality for $p=t\sqrt{n}$,
Theorem~\ref{theorem:imprPaouris} implies
Theorem~\ref{theorem:paouris} with $c=(C'+1)eC$.

Let us recall the idea underlying the proof by Paouris. Let $X\in
\R^n$ be an isotropic log-concave random vector. Let $p\sim \sqrt n$
be an integer (for example, $p = [\sqrt n]$).
Let $Y=PX$ where $P$ is an orthogonal projection of rank $p$ and let
$G\in\text {Im} P $ be a standard Gaussian vector.  By rotation
invariance, $\E |Y|^p\sim \E|\langle G/\sqrt p, Y\rangle|^{p}.$ If the
linear forms $\langle z, X\rangle$ with $ |z|=1$ had a sub-Gaussian
tail behaviour, the proof would be straightforward. But in general they
only obey a sub-exponential tail behaviour. The first step of the proof
consists of showing that there exists some $z$ for which $(\E|\langle
z, Y\rangle|^{p})^{1/p}$ is in fact small, compared to $\E |Y|$.  The
second step uses a concentration principle to show that $(\E_X|\langle
z, PX\rangle|^{p})^{1/p}$ is essentially constant on the sphere, for a
random orthogonal projection of rank $p\sim \sqrt n$, and thus
comparable to the minimum. Thus for these {\em good} projections, one
has a good estimate of $(\E |Y|^p)^{1/p}$ and the result follows by
averaging over $P$. Our proof follows the same scheme, at least for
the first step, but whereas the proof of the first step in \cite{Pa}
is the most technical part, our argument is very simple.  Then the
estimate for $\min_{|z|=1} \E|\langle z, Y\rangle|^{p}$ brings us to a
minimax problem precisely in the form answered by Gordon's inequality
(\cite{Go}).

{}Finally we would like to note that our proof can be generalized to the
case of convex measures in the sense of \cite{Bo1, Bo2}.  Of course
the proof is longer and more technical. We provide the details in
\cite{seven}.

\section{Proof of Theorem \ref{theorem:imprPaouris}}

\medskip

{}First let us notice that it is enough to prove
Theorem~\ref{theorem:imprPaouris} for symmetric log-concave random
vectors. Indeed, let $X$ be a log-concave random vector and let $X'$
be an independent copy. By Jensen's inequality we have for all $p \ge
1$,
\begin{displaymath}
(\E|X|^p)^{1/p} \le (\E|X - \E X|^p)^{1/p}+|\E X| \le
(\E|X -  X'|^p)^{1/p} +\E |X|.
\end{displaymath}
On the other hand $\E |X - X'|\le 2\E|X|$ and for $p \ge 1$ one has
$\sigma_p(X-X')\leq 2 \sigma_p(X)$. Since $X - X'$ is log-concave (see
\cite{DKH}) and symmetric, we obtain that the symmetric case proved
with a constant $C'$ implies the non-symmetric case with the constant
$C=2C' +1$.

\medskip

\begin{lemma}
\label{lemma:minimum}
Let $Y\in \R^q$ be a random vector.  Let $\|\cdot\|$ be a norm on
$\R^q$.  Then for all $p >0$,
\[
\min_{|z|=1}
(\E|\langle z, Y\rangle|^{p})^{1/p } \leq
\frac{(\E\|Y\|^{p})^{1/p}}{\E\|Y\|}\,\E|Y|.
\]
\end{lemma}

\noindent{\bf Proof:} Let  $r$ be the largest number such that
$r\|t\|\ \leq |t|$ for all  $t\in \R^q$.  Using duality, pick
$z\in\R^q$ such that
$|z|=1$ and $\|z\|_*\leq r$ (the dual norm of $\|\cdot\|$).
Then $|\langle z, t \rangle| \leq r\|t\| \leq |t|$ for all  $t\in \R^q$.
Therefore, $(\E|\langle z, Y \rangle|^{p})^{1/p} \leq r(\E\|Y\|^{p})^{1/p}$
for any $p>0$,
and the proof follows from $r \E\| Y\| \leq \E|Y|$.
\qed

\begin{lemma}
\label{lemma:ratio logconcave}
Let $Y$ be a full-dimensional symmetric log-concave $\R^q$-valued
random vector.  Then there exists a norm $\| \cdot \|$ on $\R^q$ such
that
 \[
(\E\|Y\|^{q})^{1/q} \leq 500\, \E\|Y\| .
\]
\end{lemma}

\noindent{\bf Remark.} In fact the constant 500 can be significantly
improved. To keep the presentation short and transparent we omit
the details.

\medskip

\noindent{\bf Proof:} From Borell's characterization $Y$ has an even
log-concave density $g_Y$. Thus $g_Y(0)$ is the maximum of $g_Y$.
Define a symmetric convex set by
\[
K=\{ t \in \R^q: g_Y(t) \geq 25^{-q} g_Y(0)\}.
\]
Since clearly $K$ has a non-empty interior, it is the unit ball of a
norm which we denote by $\| \cdot\|$.  On one hand, $ 1 \geq \PP(Y \in
K)=\int_{K}g_Y \geq 25^{-q}g_Y(0)\text{vol}(K), $ thus
\[
\PP(\| Y\|\leq 1/50)=\int_{K/50}g_Y \leq g_Y(0)50^{-q}\text{vol}(K)
\leq 2^{-q} \leq 1/2.
\]
Therefore
$\E\|Y\| \geq \PP(\| Y\|>1/50)/50 \geq 1/100. $
On the other hand, by the log-concavity of $g_Y$,
\[
\forall t \in \R^{q} \setminus K \qquad g_{2Y}(t)=2^{-q}g_{Y}(t/2)
\geq 2^{-q}g_{Y}(t)^{1/2}g_{Y}(0)^{1/2} \geq (5/2)^{q}g_{Y}(t).
\]
Therefore
\[
\E\| Y\|^q \leq 1+\E(\| Y\|^{q}1_{Y\in \R^{q}  \setminus K}) \leq
1+(2/5)^{q}\E\|2Y\|^{q}= 1+(4/5)^q\E\| Y\|^{q}.
\]
We conclude that  $(\E\|Y\|^{q})^{1/q}\leq  5$ and
$(\E\|Y\|^{q})^{1/q}/\E\|Y\|\leq 500.$
\qed

\begin{lemma}
\label{lemma:minimax}
Let $n, q\geq 1$ be integers and $p\geq 1$. Let $X$ be an
$n$-dimensional random vector, $G$ be a standard Gaussian vector in
$\R^n$ and $\Gamma$ be an $n\times q$ standard Gaussian matrix. Then
\[
  (\E |X|^p)^{1/p}\leq \alpha_p^{-1}\left(\E\min_{|t|=1}|||\Gamma
  t|||+(\alpha_p+ \sqrt q)\,\sigma_p(X) \right),
\]
where $|||z|||=(\E|\langle z, X\rangle|^{p})^{1/p}$ and $\alpha_p^p$
is the $p$-th moment of an $N(0,1)$ Gaussian random variable
(so that
$\lim_{p\to\infty}\big(\alpha _p/\sqrt p\big) = 1/\sqrt e$).
\end{lemma}

\noindent {\bf Proof:}
 By rotation invariance,
$\E|\langle G, X\rangle|^{p}=\alpha_p^p \,\E |X|^p$.
Notice that
$$
  \sigma^2 :=  \sup _{|||t|||_* \leq 1} \E |\langle G, t \rangle|^2 =
   \sup _{|||t|||_* \leq 1} |t|^2 = \sigma_p^2 (X),
$$
where $|||\cdot |||_*$ denotes the norm on $\R^n$ dual to the norm
$|||\cdot |||$. Denote the median of $|||G|||$ by $M_G$.
The classical deviation inequality for a norm of a Gaussian vector
(\cite{Bo3}, \cite{SC}, see also \cite{L}, Theorem 12.2)  states
\[
  \forall s \geq 0\qquad \PP \left( \big| \, ||| G||| - M_G \big|
  \geq s \right) \leq 2\, \int _{s/\sigma}^{\infty} \exp \big({-t^2/2}\big)
 \ \frac{dt}{\sqrt{2\pi}}
\]
and since $M_G \leq \E||| G|||$  (\cite{Kw}, see also \cite{L}, Lemma~12.2)
this implies
\[
  (\E |X|^p)^{1/p}= \alpha_p^{-1} (\E||| G|||^p)^{1/p} \leq
  \alpha_p^{-1}\big(\E||| G||| + \alpha _p \sigma_p(X) \big)
\]
(cf. \cite{LMS}, Statement 3.1).

The  Gordon minimax lower bound (see \cite{Go},  Theorem 2.5)
states that for any norm $||| \cdot  |||$
\[
  \E ||| G|||\leq \E\min_{|t|=1}|||\Gamma t|||+ \left(\E |H|\right)
  \, \max_{|z|=1}|||z||| \leq \E\min_{|t|=1}|||\Gamma t||| +
  \sqrt{q} \, \sigma_p(X),
\]
where  $H$ is a  standard Gaussian vector in $\R^q$.
This concludes the proof.
\qed

\medskip

\noindent{\bf Proof of Theorem~\ref{theorem:imprPaouris}:} Assume
that $X$ is log-concave symmetric. We use the notation of Lemma
\ref{lemma:minimax} with $q$ the integer such that $p\leq q< p+1$. We
first condition on $\Gamma$.  Let $Y=\Gamma^*X$. Note that $Y$ is
log-concave symmetric and that
\[
  |||\Gamma t|||=(\E_X|\langle \Gamma t, X\rangle|^{p})^{1/p}=
  (\E_X|\langle t, \Gamma^*X\rangle|^{p})^{1/p} .
\]
If $\Gamma^*X$ is supported by a hyperplane then
$\min_{|t|=1}(\E_X|\langle t, \Gamma^*X\rangle|^{p})^{1/p}=0$.
Otherwise Lemma \ref{lemma:ratio logconcave} applies and combined with
Lemma \ref{lemma:minimum} gives that
\[
  \min_{|t|=1} |||\Gamma t||| \leq \min_{|t|=1}(\E_X|\langle t
  ,\Gamma^*X\rangle|^{p})^{1/p}\leq  500\,\E_X|\Gamma^*X|.
\]
By
taking  expectation over $\Gamma$
we get
\[
\E\min_{|t|=1} |||\Gamma t |||\leq
500 \,\E|\Gamma^*X|=500\, \E |H| \,\E|X|
\leq 500 \, \sqrt q \, \E|X|,
\]
where $H\in\R^q$ is  a  standard Gaussian vector.
Applying Lemma~\ref{lemma:minimax} we obtain
$$
  (\E |X|^p)^{1/p} \leq 500\, \alpha_p^{-1}\, \sqrt q\, \E|X|
  +   (1+\alpha_p^{-1}\sqrt q) \sigma_p(X) .
$$
 This implies the desired result, since $q\le p+1$ and hence
$\alpha ^{-1}_p \sqrt{q} \le c $ for some numerical constant $c$
 (recall that
$\lim_{p\to\infty}\big(\alpha _p/\sqrt p\big) = 1/\sqrt e$).
\qed

\address


\begin{thebibliography}{99}

\bibitem{seven} R.~Adamczak, O.~Gu\'{e}don, R.~Lata{\l}a, A.~E.~Litvak,
  K.~Oleszkiewicz, A.~Pajor and N.~Tomczak-Jaegermann, \emph{Moment
    estimates for convex measures}, preprint.

\bibitem{ALLPT} R.~Adamczak, R.~Lata{\l}a, A.~E.~Litvak, A.~Pajor and
  N.~Tomczak-Jaegermann, \emph{Tail estimates for norms of sums of
    log-concave random vectors}, preprint,
  http://arxiv.org/abs/1107.4070.

\bibitem{cras_allpt} R.~Adamczak, R.~Lata{\l}a, A.~E.~Litvak, A.~Pajor
  and N.~Tomczak-Jaegermann, \emph{Geometry of log-concave Ensembles
    of random matrices and approximate reconstruction},
  C.R. Math. Acad. Sci. Paris, {\bf 349} (2011), 783--786.

\bibitem{bmp} R.~E.~Barlow, A.~W.~Marshall and F.~Proschan, \emph
    {Properties of probability distributions with monotone hazard
      rate}, Ann. Math. Statist., {\bf 34} (1963), 375-–389.

\bibitem{Bo1} C.~Borell, {\em Convex measures on locally convex
    spaces}, Ark. Math., {\bf 12} (1974), 239--252.

\bibitem{Bo2} C.~Borell, {\em Convex set functions in d-space},
  Period. Math. Hungar., {\bf 6} (1975), 111--136.

\bibitem{Bo3} C.~Borell, {\em The Brunn-Minkowski inequality in Gauss
    space}, Invent. Math., {\bf 30} (1975), 207--216.

\bibitem{DKH} Ju.~S.~Davidovic, B.~I.~Korenbljum and B.~I.~Hacet, {\em
    A certain property of logarithmically concave functions}, Soviet
  Math. Dokl., {\bf 10} (1969), 447--480; translation from
  Dokl. Akad. Nauk SSSR {\bf 185} (1969), 1215--1218.

\bibitem{Go} Y.~Gordon, {\em Some inequalities for Gaussian processes
    and applications}, Israel J. Math., {\bf 50} (1985), 265--289.

\bibitem{Kw} S.~Kwapie{\'n}, {\em A remark on the median and the
    expectation of convex functions of Gaussian vectors}, Probability
  in Banach spaces, 9 (Sandjberg, 1993), 271--272, Progr. Probab., 35,
  Birkhäuser Boston, Boston, MA, 1994.

\bibitem{klo} S.~Kwapie{\'n}, R.~Lata{\l}a and K.~Oleszkiewicz,
  \emph{Comparison of moments of sums of independent random variables
    and differential inequalities}, J. Funct. Anal., {\bf 136} (1996),
  258-­268.

\bibitem{L} M.~A.~Lifshits, {\em Gaussian random functions.}
  Mathematics and its Applications; 322. Kluwer Academic Publishers,
  Dordrecht, 1995.


\bibitem{LMS} A.~E.~Litvak, V.~D.~Milman and G.~Schechtman, {\em
    Averages of norms and quasi-norms,} Math. Ann., {\bf 312} (1998),
  95--124.

\bibitem{Pa} G.~Paouris, {\em Concentration of mass on convex bodies},
  Geom. Funct. Anal., {\bf 16} (2006), 1021--1049.

\bibitem{SC} V.~N.~Sudakov, B.~S.~Cirel'son, {\em Extremal properties of
    half-spaces for spherically invariant measures},
  J. Sov. Math. {\bf 9} (1978), 9--18; translation from
  Zap. Nauchn. Sem. Leningrad. Otdel. Mat. Inst. Steklov (LOMI) {\bf
    41} (1974), 14--24.


\end{thebibliography}
\end{document}